\documentclass{article}
     \usepackage[english]{babel}
     \usepackage{hyperref, amsmath, amsthm, amssymb, amsfonts, xcolor, mathdots, breqn, authblk}

\newcommand{\bit}{\begin{itemize}}
\newcommand{\eit}{\end{itemize}}

\newcommand{\bmat}{\begin{dmath}}
\newcommand{\emat}{\end{dmath}}

\newcommand{\bmatn}{\begin{dmath*}}
\newcommand{\ematn}{\end{dmath*}}

\newtheorem{thm}{Theorem}
\newtheorem{prop}[thm]{Proposition}
\newtheorem{ex}{Example}

\DeclareMathOperator{\ric}{Ric}
\DeclareMathOperator{\e}{e}

\newcommand*{\defeq}{\mathrel{\vcenter{\baselineskip0.5ex \lineskiplimit0pt
                     \hbox{\scriptsize.}\hbox{\scriptsize.}}}%
                     =}

\title{Ricci flat Finsler metrics by warped product}

\author{Patricia Marcal}
\author{Zhongmin Shen}

\affil[]{Department of Mathematical Sciences \\ IUPUI, USA}

\affil[ ]{\textit {\{pmarcal, zshen\}@iupui.edu}}

\date{\vspace{-5ex}}

\begin{document}

\maketitle

\bibliographystyle{abbrv}

\begin{abstract}
In this work, we consider a class of Finsler metrics using the warped product notion introduced by Chen, S. and Zhao in \cite{chen:shen:zhao}, with another ``warping'', one that is consistent with static spacetimes. We will give the PDE characterization for the proposed metrics to be Ricci-flat and explicitly construct two non-Riemannian examples.

\end{abstract}

\section{Introduction}

If $(M,d s_1^2)$, $(N,d s_2^2)$ are Riemannian manifolds, then a warped product is the manifold $M \times N$ endowed with a Riemannian metric of the form
\begin{equation}
\label{warped prod riem}
    d s^2 = d s_1^2 + f^2 d s_2^2 \, ,
\end{equation}
where $f$ is a smooth function depending on the coordinates of $M$ only; said a warping function. This notion, called by \emph{warped product}, must be credited to Bishop and O'Neill \cite{bishop:oneill}. However, years earlier, metrics in the form of (\ref{warped prod riem}) were being studied with different names; in \cite{kruch}, for instance, they were called semi-reducible Riemannian spaces.

The class of warped product manifolds has shown itself to be rich, both wide and diverse, playing important roles in differential geometry as well as in physics. To illustrate, Bishop and O'Neill introduced warped products in \cite{bishop:oneill} as means to construct a large class of complete Riemannian manifolds with negative curvature. For this reason, it seems valuable to study notions of warped product metrics without the quadratic restriction, in the setting of Finsler geometry. Notably, progress in this direction has been stimulated by efforts to expand general relativity, such as the work of Asanov (e.g. \cite{asanov85}, \cite{asanov91}, \cite{asanov98}), which later motivated Kozma, Peter and Varga to study product manifolds $M \times N$ endowed with a Finsler metric
\begin{equation}
\label{warped prod 1}
    F = \sqrt{F_1^{2} + f^{2}F_{2}^{2}} \, ,
\end{equation}
called warped product, where $(M,F_1)$, $(N,F_2)$ are Finsler manifolds and $f$ is a smooth function on $M$ (see \cite{kozma:peter:varga}). Following the definition of Beem \cite{beem}, one may take $L = F^2$ to consider pseudo-Finsler metrics. For example, if $(M, F_1)$ is a $3$-dimensional Finsler manifold and $(\mathbb{R},F_2)$ is a Minkowski space, then
\begin{equation}
\label{static spacetime 1}
    L = f^2 F_2^2 - F_1^2
\end{equation}
is a Finsler metric with Lorentz signature, and $(\mathbb{R}\times M, L)$ may be regarded as a Finsler static spacetime. This is the case for \cite{li:chang}, where Li and Chang studied metrics in the form of (\ref{static spacetime 1}), given on coordinates $((t,r,\theta , \varphi),(y^t,y^r,y^\theta , y^\varphi))$ of the tangent bundle by
\begin{equation*}
    L = f^2 (y^t)^2 - \left[g^2(y^r)^2 + r^2\overline{F}^2\right] \, ,
\end{equation*}
with $\overline{F}$ a Finsler metric on coordinates $(\theta, \varphi,y^\theta , y^\varphi)$ and $f,g$ functions of $r$. They suggested the vacuum field equation for Finsler spacetime is equivalent to the vanishing of the Ricci scalar, and obtained a non-Riemannian exact solution similar to the Schwarzschild metric. Their results display viability to the explict construction of Ricci-flat Finsler metrics by warped product.

Recently, Chen, S. and Zhao have considered product manifolds $\mathbb{R}\times M$ with Finsler metrics arising from warped products in the following way: if $(M,\alpha^2)$, $(\mathbb{R},d t^2)$ are Riemannian manifolds, then $F^2 = d t^2 + f^2(t)\alpha^2$ is a warped product, which may be rewritten as $F = \alpha\sqrt{\left(\frac{d t}{\alpha}\right)^2 + f^2(t)}$. Letting $z = \frac{d t}{\alpha}$, they defined a class of Finsler metrics by
\begin{equation}
\label{warped prod 2}
    F = \alpha \sqrt{\phi(z,t)} \, ,
\end{equation}
also called warped product, where $\phi$ is a suitable function on $\mathbb{R}^2$, (see \cite{chen:shen:zhao}). 

In the present work, we wish to consider Finsler metrics of similar type as (\ref{warped prod 2}), with another ``warping'', one that is consistent with the form of metrics modeling (Riemannian) static spacetimes and simplified by spherical symmetry over spatial coordinates, which emerged from the Schwarzschild metric in isotropic rectangular coordinates $(t,x^1,x^2,x^3)$:
\begin{equation}
\label{sch iso coord}
    d s^2 = \frac{\left(1 - \frac{m}{4\rho}\right)^2}{\left(1 + \frac{m}{4\rho}\right)^2} c^2 d  t^2 - \left(1 + \frac{m}{4\rho}\right)^4 \left[(d x^1)^2 + (d  x^2)^2 + (d x^3)^2\right] \, ,
\end{equation}
where $\rho=\sqrt{(x^1)^2+(x^2)^2+(x^3)^2}$ (see for example \cite{edd}, p. 93). Letting $z = \frac{d t}{\alpha}$ and $\alpha = \sqrt{(d x^1)^2 + (d x^2)^2 + (d x^3)^2}$, the Schwarzschild metric (\ref{sch iso coord}) is written as
$$d s = \alpha \sqrt{ \frac{\left(1 - \frac{m}{4\rho}\right)^2}{\left(1 + \frac{m}{4\rho}\right)^2} c^2 z^2 - \left(1 + \frac{m}{4\rho}\right)^4 } \, .$$

For $\mathbb{R}$, $\mathbb{R}^n$ with their Euclidean metrics $d t^2$, $\alpha^2$ (respectively), define a (positive-definite) Finsler metric on $\mathbb{R}\times\mathbb{R}^n$ by
\begin{equation*}
    F = \alpha \sqrt{\phi(z,\rho)} \, ,
\end{equation*}
where $z = \frac{d t}{\alpha}$, $\rho = \vert \overline{x}\vert$ for $\overline{x}\in \mathbb{R}^n$, and $\phi$ a suitable function on $\mathbb{R}^2$. We give the PDE characterization for the proposed metrics to be Ricci-flat:
\begin{thm}
\label{thm}
For $n\geq 2$, $F = \alpha\sqrt{\phi(z,\rho)}$ is Ricci-flat if and only if $P(z,\rho) = Q(z,\rho) = 0$. Furthermore, the Ricci-flat condition is weaker when $n=1$; namely, $P(z,\rho) + \rho^2Q(z,\rho) = 0$.
\end{thm}

Here $P$, $Q$ are functions of $\phi$ and its derivatives, described by equations (\ref{ricci comp}).

Next, we construct three examples. The first presents the Riemannian case to corroborate computations, the second are $m$-th root metrics, and the third consists of Randers norm. For $n\geq 3$, the non-Riemannian solutions are:
\begin{align*}
    \phi(z, \rho) &= (Az^m + B\rho^{-2m})^{\frac{2}{m}} \, , \; A,B > 0 &(\ref{m3sol})\\
    \phi(z, \rho) &= (\sqrt{A z^2 + B\rho^{-4}} + \varepsilon\sqrt{A}z)^2 \, , \; A,B > 0 \, , \; 0 < \vert \varepsilon \vert < 1 \, , \; C\in\mathbb{R}  &(\ref{r3sol})
\end{align*}
Whenever possible, we describe solutions with Lorentz signature, for which the $4$-dimensional metrics may also be studied as Finsler static spacetimes satisfying the vacuum field equation proposed in \cite{li:chang}.

\section{Geometric Quantities}

Set $M = \mathbb{R}\times\mathbb{R}^n$ with coordinates on $TM$
\begin{align*}
    x=(x^0,\overline{x}) &, \; \overline{x} = (x^1,\ldots , x^n) \, , \\
    y=(y^0,\overline{y}) &, \; \overline{y} = (y^1,\ldots , y^n) \, ;
\end{align*}
and consider a Finsler metric
\begin{equation}
\label{def}
    F = \alpha\sqrt{\phi(z,\rho)}\, ,
\end{equation}
where $\alpha= \vert\overline{y}\vert$, $z=\frac{y^0}{\vert\overline{y}\vert}$ and $\rho=\vert\overline{x}\vert$. Throughout our work, the following convention for indices is adopted: A, B, ... range from $0$ to $n$; i, j, ... range from $1$ to $n$.

This construction is the same as \cite{chen:shen:zhao} but for the ``warping''. Consequently, any calculations involving $F$ and its derivatives of any degree with respect to $y^A$ only will be similar in form to the calculations in \cite{chen:shen:zhao}, e.g. the fundamental form. The effects of the warping only appear when derivatives of $F$ with respect to $x^A$ are involved, e.g. spray coefficients. So the Hessian matrix, $g_{AB} = \frac{1}{2}[F^2]_{y^A y^B}$, is
\begin{equation}
\label{hessian matrix}
    (g_{AB}) = \left(
    \begin{array}{c|c}
    \frac{1}{2}\phi_{zz} & \frac{1}{2}\Omega_z\frac{y^j}{\alpha} \\
    \hline
    \frac{1}{2}\Omega_z\frac{y^i}{\alpha} & \frac{1}{2}\Omega\delta_{ij} - \frac{1}{2}z\Omega_z\frac{y^i y^j}{\alpha^2}
    \end{array}
    \right) \, ,
\end{equation}
where
\begin{equation}
    \Omega := 2\phi -z\phi_z\, ,
\end{equation}
and the same argument as \cite{chen:shen:zhao} to verify non-degeneracy of $F$ applies. It actually simplifies, because $\alpha$ is the Euclidean metric here. Thus,
$$\det (g_{AB}) = \frac{1}{2^{n+1}}\Omega^{n-1}\Lambda \, ,$$
where
\begin{equation}
    \Lambda:=\phi_{zz}(\Omega - z\Omega_z) - \Omega_z^2 = 2\phi\phi_{zz} - \phi_z^2\, ,
\end{equation}
and:

\begin{prop}[Prop.4.1, \cite{chen:shen:zhao}]
$F=\alpha\sqrt{\phi(z,\rho)}$ is strongly convex if and only if $\Omega, \Lambda > 0$.
\end{prop}

Moreover, letting $L = \alpha^2\phi(z,\rho)$, one has  a metric with Lorentz signature $(+,-,\ldots,-)$ if $\Omega,\Lambda < 0$, or $(-,+,\ldots,+)$ if $\Omega > 0$ and $\Lambda < 0$.

Henceforth, assume $F$ is non-degenerate. In this case, the inverse of $(g_{AB})$ is
\begin{equation}
\label{inverse matrix}
    (g^{AB}) = \left(
    \begin{array}{c|c}
    \frac{2}{\Lambda}(\Omega - z\Omega_z) & -\frac{2}{\Lambda}\Omega_z\frac{y^j}{\alpha} \\
    \hline
    -\frac{2}{\Lambda}\Omega_z\frac{y^i}{\alpha} & \frac{2}{\Omega}\delta^{ij} + \frac{2\phi_z\Omega_z}{\Omega\Lambda}\frac{y^i y^j}{\alpha^2}
    \end{array}
    \right) \, .
\end{equation}

The spray coefficients $G^C = \frac{1}{4}g^{CA}\left([F^2]_{y^A x^B}y^B - [F^2]_{x^A} \right)$ are:
\begin{subequations}
\begin{align}
    G^0 &= (U + z V)(x^m y^m)\alpha \, ,    \label{g0} \\
    G^i &= (V + W) y^i(x^m y^m) - W x^i\alpha^2 \, ,    \label{gi}
\end{align}
\end{subequations}
where
\begin{subequations}
\begin{align}
    U &:= \frac{1}{2\rho\Lambda}(2\phi\phi_{z\rho} - \phi_z\phi_{\rho})\, ,    \label{u}\\
    V &:= \frac{1}{2\rho\Lambda}(\phi_{\rho}\phi_{zz} - \phi_z\phi_{z\rho})\, ,    \label{v}\\
    W &:= \frac{1}{2\rho\Omega}\phi_{\rho}\, .    \label{w}
\end{align}
\end{subequations}

The flag curvature tensor by Berwald's formula $$K_B^C = 2[G^C]_{x^B} - [G^C]_{x^A y^B} y^A + 2G^A[G^C]_{y^A y^B} - [G^C]_{y^A}[G^A]_{y^B}$$ gives
\begin{dgroup}[frame={0pt},framesep={4pt}]
\bmat
K_0^0 = \left[ \rho^2(U + z V)W_z - (2\rho^2W +1)(U_z + V + z V_z) \right] \alpha^2 + \left[ 2(V + W)(U_z + V + z V_z) - (V_z + W_z)(U + z V) + 2U(U_{zz} + 2V_z + z V_{zz}) - \frac{1}{\rho}(U_{z\rho} + V_{\rho} + z V_{z\rho}) - (U_z + V + z V_z)^2 - (U  - z U_z - z^2V_z)V_z \right](x^m y^m)^2
\emat

\bmat
K_j^i = - \left[ 2W + (2\rho^2W + 1)(V + W) \right]\alpha^2\delta^i_j +\left[ (V + W)^2 + 2U(V_z + W_z) - \frac{1}{\rho}(V_{\rho} + W_{\rho}) \right](x^m y^m)^2\delta^i_j + \left[ 2W(2W -z W_z) + W_z(U - zW) -\frac{2}{\rho}W_{\rho} \right]\alpha^2x^i x^j + \left[ (V + W) + z(V_z + W_z)(2\rho^2W + 1) + (\rho^2(V + W) + 1)(2W - z W_z) \right]y^i y^j - \left[ 2z U (V_{zz} + W_{zz}) + (3U - z U_z - z V + 5z W)(V_z + W_z) - \frac{z}{\rho}(V_{z\rho} + W_{z\rho}) \right] (x^m y^m)^2\frac{y^i y^j}{\alpha^2} + \left[ - (2W - z W_z)^2 -2U(W_z - z W_{zz}) + \frac{1}{\rho}(2W_{\rho} - z W_{z\rho}) + W_z(U - z U_z + z^2 W_z) \right](x^m y^m) x^i y^j + \left[ -(V + W)^2 + (V_z + W_z)(U + 3z W) + \frac{1}{\rho}(V_{\rho} + W_{\rho}) \right](x^m y^m) x^j y^i
\emat

\bmat
 K_j^0 = z\left[ (2\rho^2 W + 1)(V + U_z + z V_z) - \rho^2 W_z(U + z V) \right]\alpha y^j  + \left[ z(U + z V)(V_z + W_z) - 2z U(U_{zz} + 2V_z + z V_{zz}) + (U - z U_z - z^2 V_z)(5W - U_z) - \frac{1}{\rho}(U_{\rho} - z U_{z\rho} - z^2 V_{z\rho}) \right](x^m y^m)^2\frac{y^j}{\alpha} + \left[ (U + z V)(U_z - V + z V_z - 2W) + (V - 3W)(U - z U_z - z^2 V_z) + \frac{1}{\rho}(U_{\rho} + z V_{\rho}) \right](x^m y^m)\alpha x^j
\emat

\bmat
 K_0^i = \left[ \rho^2 W_z(V - W) - (2\rho^2W + 1)V_z \right]\alpha y^i + \left[ (2W - V - U_z)(V_z + W_z) + 2U(V_{zz} + W_{zz}) -\frac{1}{\rho}(V_{z\rho} + W_{z\rho}) \right](x^m y^m)^2 \frac{y^i}{\alpha} + \left[ (U_z - W)W_z - 2U W_{zz} + \frac{1}{\rho}W_{z\rho} \right](x^m y^m)\alpha x^i
\emat
\end{dgroup}

After simplification, the Ricci curvature is:
\bmatn
\ric = \sum K^A_A \\
        = \left[ - (2\rho^2 W + 1)(U_z + n V + (n-3) W) - 2(n W + \rho W_\rho - \rho^2 W_z(U - z W)) \right]\alpha^2 + \left[ 2U(U_{zz} + n V_z + (n-2) W_z) -\frac{1}{\rho}(U_{z\rho} + n V_{\rho} + (n-3) W_{\rho}) + n V(V + 2W) + W((n-5) W + 2z W_z) + U_z(2W - U_z) \right](x^m y^m)^2
\ematn

Let the Ricci curvature components be
\begin{dgroup}[frame={0pt},framesep={4pt}]
\label{ricci comp}
\bmat
 P(z,\rho) \defeq - (2\rho^2 W + 1)(U_z + n V + (n-3) W) - 2(n W + \rho W_\rho - \rho^2 W_z(U - z W))
\emat
\bmat
 Q(z, \rho) \defeq 2U(U_{zz} + n V_z + (n-2) W_z) -\frac{1}{\rho}(U_{z\rho} + n V_{\rho} + (n-3) W_{\rho}) + n V(V + 2W) + W((n-5) W + 2z W_z) + U_z(2W - U_z)
\emat
\end{dgroup}

So
\bmat
\label{ricci}
\ric = P\left(\frac{y^0}{\vert\overline{y}\vert},\vert\overline{x}\vert\right)\langle \overline{y} , \overline{y} \rangle + Q\left(\frac{y^0}{\vert\overline{y}\vert},\vert\overline{x}\vert\right)\langle \overline{x} , \overline{y} \rangle^2 \\
        = \left\langle P\left(\frac{y^0}{\vert\overline{y}\vert},\vert\overline{x}\vert\right)\overline{y} + Q\left(\frac{y^0}{\vert\overline{y}\vert},\vert\overline{x}\vert\right)\langle\overline{x},\overline{y}\rangle\overline{x} , \overline{y} \right\rangle
\emat

\begin{proof}[{\bf Proof of Theorem \ref{thm}}]
Suppose $\ric = 0$. Let $e_i$ denote the $n$-dimensional vector with $1$ in the $i^{th}$ entry and zeros elsewhere. Take $\overline{y} = e_i$ and $\overline{x} = \rho e_j$ for $\rho \geq 0$. By equation (\ref{ricci}),
$$P\left(y^0,\rho\right) + Q\left(y^0,\rho\right)\rho^2\delta^{ij} = 0 \, , \; \forall i,j \, .$$
For $n \geq 2$, pick $i\neq j$ to get $P(y^0,\rho) = 0$. Now set $i=j$ to conclude $Q(y^0,\rho) = 0$ for $\rho \neq 0$. Finally, $Q(y^0,0) = 0$ by continuity. The remaining assertions are clear.
\end{proof}

The above proof suggests metrics $F$ that are singular on $(x^0,0)$ or  metrics $F$ defined on $\mathbb{R}\times \mathbb{R}^n\setminus\{0\}$ should also be considered. This becomes evident on the examples bellow.

\section{Examples}

\begin{ex}[Riemannian metrics]
Suppose $\phi(z, \rho) = \e^{f(\rho)}z^2 + \e^{g(\rho)}$.

So $\Omega = 2\e^{g}$, $\Lambda = 4\e^{f + g}$ and $F = \alpha\sqrt{\phi}$ gives a positive-definite Riemannian metric.

The Ricci curvature components are:
\begin{align*}
    P &= -\frac{1}{4\rho}\left[ p_2\e^{f-g}z^2  + p_0 \right] \\
    Q &= -\frac{1}{4\rho^3}q_0
\end{align*}
where
\begin{dgroup*}
\bmatn
p_2 = 2\rho f^{\prime\prime} + \rho(f^{\prime})^2 + (n-2)\rho f^{\prime}g^{\prime} + 2(n-1)f^{\prime}
\ematn
\bmatn
p_0 = 2\rho g^{\prime\prime} + (n-2)\rho(g^{\prime})^2 + \rho f^{\prime}g^{\prime} + 2f^{\prime} + 2(2n-3)g^{\prime}
\ematn
\bmatn
q_0 = 2\rho f^{\prime\prime} + 2(n-2)\rho g^{\prime\prime} + \rho(f^{\prime})^2 - 2\rho f^{\prime}g^{\prime} - (n-2)\rho (g^{\prime})^2 - 2f^{\prime} - 2(n-2)g^{\prime}
\ematn
\end{dgroup*}

By independence of $z$ and $\rho$, the Ricci-flat equations for $n\geq 2$ become $p_2 = p_0 = q_0 = 0$. Taking $q_0 - p_2 + n p_0 = 0$ yields:
$$4(n-1)\rho g^{\prime\prime} + (n-2)(n-1)\rho(g^{\prime})^2 + 4(n-1)^2g^{\prime} = 0 $$

For $n\geq 3 :$
\begin{equation}
\label{Riemann edo}
    4\rho g^{\prime\prime} + (n-2)\rho(g^{\prime})^2 + 4(n-1)g^{\prime} = 0 
\end{equation}
If $g^{\prime} = 0$, then (\ref{Riemann edo}) is trivially satisfied. So $g(\rho)=\overline{B}$ constant is a solution. Otherwise, (\ref{Riemann edo}) is a Bernoulli differential equation in $g^{\prime}$, which can be transformed to a linear ODE by letting $u:= (g^{\prime})^{-1}$. The equation reduces to:
$$4\rho u^{\prime} - (n-2)\rho - 4(n-1)u = 0$$
Its solution gives $g(\rho) = \ln(B\vert \rho^{2-n} + C\vert^{\frac{4}{n-2}})$, for $B,C\in \mathbb{R}$ constants with $B>0$.

To find $f$, substitute $g$ in $p_0 = 0$. When $g(\rho)=\overline{B}$, $f^{\prime} = 0$ and $f$ is constant also, say $f(\rho) = \overline{A}$. For $g(\rho) = \ln(B\vert \rho^{2-n} + C\vert^{\frac{4}{n-2}})$, equation $p_0= 0$ gives:
$$f^{\prime} = \frac{4(n-2)C\rho^{1-n}}{(C-\rho^{2-n})(C+\rho^{2-n})}$$
So $f(\rho) = \ln \left[A\left(\frac{C-\rho^{2-n}}{C+\rho^{2-n}}\right)^2\right]$, for some constant $A>0$.

Therefore, when $n\geq 3$, solutions are:
\begin{subequations}
\begin{align}
    \phi(z, \rho) &= Az^2 + B \, , \; A,B > 0 \\
    \phi(z, \rho) &= A\left(\frac{C-\rho^{2-n}}{C+\rho^{2-n}}\right)^2 z^2 + B\vert \rho^{2-n} + C\vert^{\frac{4}{n-2}} \, , \; A,B > 0 , \, C\in\mathbb{R}
\end{align}
\end{subequations}

For $n = 2$, equation (\ref{Riemann edo}) still holds, but it is already linear:
$$\rho g^{\prime\prime} + g^{\prime} = 0$$
So $g(\rho) = \ln (B\vert \rho\vert^C)$, for $B,C\in \mathbb{R}$ constants with $B>0$. Substitute $g$ in $p_0 = 0$ to get:
$$(C+2)f^{\prime} = 0$$
If $C\neq -2$, then $f^{\prime} = 0$. So $f(\rho) = \overline{A}$. When $C = -2$, equation $p_2 = 0$ yields:
$$2\rho f^{\prime\prime} + \rho(f^{\prime})^2 + 2f^{\prime} = 0$$
If $f^{\prime} = 0$, the above equation is trivially satisfied and then $f(\rho) = \overline{A}$. Else, it is a Bernoulli equation in $f^{\prime}$. As before, let $u:= (f^{\prime})^{-1}$ to get a linear ODE:
$$2\rho u^{\prime} - \rho - 2u = 0$$
It gives $f(\rho) = \ln (A_1 + A_2\ln\vert\rho\vert)^2$, for real constants $A_1 , A_2$.

Thus, for $n=2$, solutions are:
\begin{subequations}
\begin{align}
    \phi(z, \rho) &= Az^2 + B\vert\rho\vert^C \, , \; A,B > 0 , \, C\in\mathbb{R}\setminus\{-2\} \\
    \phi(z, \rho) &= (A_1 + A_2\ln\vert\rho\vert)^2z^2 + B\rho^{-2} \, , \; A_1,A_2\in\mathbb{R}, \, B > 0 
\end{align}
\end{subequations}

For $n = 1$, the Ricci-flat condition gives $p_2 = p_0 + q_0 = 0$, by independence of $z$ and $\rho$. This gives:
$$2 f^{\prime\prime} + (f^{\prime})^2 - f^{\prime}g^{\prime} = 0$$
So either $f(\rho) = \overline{A}$ and $g$ is an arbitrary smooth function of $\rho$, or $g = \ln(f^{\prime})^2 + f + \overline{B}$ for any smooth function $f$ of $\rho$.

Hence, solutions for $n=1$ are:
\begin{subequations}
\begin{align}
    \phi(z, \rho) &= Az^2 + \e^{g(\rho)} \, , \; A > 0 \, , \; g\in C^{\infty} \\
    \phi(z, \rho) &= \e^{f(\rho)}\left(z^2 + B[f^{\prime}(\rho)]^2 \right) , \; f\in C^{\infty} , \, B > 0
\end{align}
\end{subequations}

Finally, if $\phi(z, \rho) = \e^{f(\rho)}z^2 - \e^{g(\rho)}$, then $\Omega = - 2\e^{g}$ and $\Lambda = - 4\e^{f + g}$. So the associated metric $L = \alpha^2\phi$ has Lorentz signature $(+,-,\ldots,-)$. In this case, the Ricci curvature components are:
\begin{align*}
    P &= \frac{1}{4\rho}\left[ p_2\e^{f-g}z^2  - p_0 \right] \\
    Q &= -\frac{1}{4\rho^3}q_0
\end{align*}
where $p_2$, $p_0$ and $q_0$ are as before. Thus, by independence of $z$ and $\rho$, the Ricci-flat equations reduce to the same system as the positive-definite case.
\end{ex}

\begin{ex}[$m^{th}$-root metrics]
If $\phi(z,\rho) = \left(\e^{f(\rho)}z^{m} + \e^{g(\rho)}\right)^{\frac{2}{m}}$ for an even integer $m > 2$,  then $\Omega = \frac{2\e^g}{(\e^f z^m + \e^g)^{1-\frac{2}{m}}}$ and $\Lambda = \frac{4(m-1)\e^{f+g}z^{m-2}}{(\e^f z^m + \e^g)^{2(1-\frac{2}{m})}}$. So $F = \alpha\sqrt{\phi}$ is a positive-definite $m^{th}$-root metric.

The Ricci curvature components are:
\begin{align*}
    P &= -\frac{1}{2m^2(m-1)\rho}\left[ p_{2m}\e^{2(f-g)}z^{2m} + p_m\e^{f-g}z^m  + p_0 \right] \\
    Q &= \frac{1}{4m^2(m-1)^2\rho^3}\left[ q_{2m}\e^{2(f-g)}z^{2m} - q_m\e^{f-g}z^m  - q_0 \right] 
\end{align*}
where
\begin{dgroup*}
\bmatn
p_{2m} = (m-2)(m+n-2)\rho(f^{\prime})^2
\ematn
\bmatn
p_m = 2m(m-1)\rho f^{\prime\prime} + m(m-1)\rho(f^{\prime})^2 + (n-2)(3m-4)\rho f^{\prime}g^{\prime} + m[(n-2)(3m-4)+2(m-1)]f^{\prime}
\ematn
\bmatn
p_0 = 2m(m-1)\rho g^{\prime\prime} + 2(m-1)(n-2)\rho(g^{\prime})^2 + m\rho f^{\prime}g^{\prime} + m^2f^{\prime} + 2m(m-1)(2n-3)g^{\prime}
\ematn
\bmatn
q_{2m} = (m-2)[2m^2 + (n-2)(3m-2)]\rho(f^{\prime})^2
\ematn
\bmatn
q_m = 2(m-2)[ m(m-1)(n-2)\rho f^{\prime\prime} - m(m+n-1)\rho(f^{\prime})^2 + 2(n-2)(m-1)\rho f^{\prime}g^{\prime} + m(m-1)(n-2)f^{\prime}]
\ematn
\bmatn
q_0 = 2m^2(m-1)\rho f^{\prime\prime} + 4m(m-1)^2(n-2)\rho g^{\prime\prime} + m^2\rho(f^{\prime})^2 - 4m(m-1)\rho f^{\prime}g^{\prime} - 4(m-1)^2(n-2)\rho (g^{\prime})^2 - 2m^2(m-1)f^{\prime} - 4m(m-1)^2(n-2)g^{\prime}
\ematn
\end{dgroup*}

By independence of $z$ and $\rho$, the Ricci-flat equations for $n\geq 2$ are $p_{2m} = p_m = p_0 = q_{2m} = q_m = q_0 = 0$. Since $m > 2$, $p_{2m} = q_{2m} = 0$ imply $f^{\prime} = 0$, and equations $p_m = q_m = 0$ are automatically satisfied. The remaining equations reduce to:
\begin{equation}
\label{m p eq}
    m\rho g^{\prime\prime} + (n-2)\rho (g^{\prime})^2 + (2n-3)mg^{\prime} = 0
\end{equation}
\begin{equation}
\label{m q eq}
      (n-2)[m\rho g^{\prime\prime} - \rho(g^{\prime})^2 - mg^{\prime}] = 0
\end{equation}
So $f(\rho) = \overline{A}$ and $g(\rho)$ must be determined from the above equations.

For $n \geq 3$, combine equations (\ref{m p eq}) and (\ref{m q eq}) to eliminate $g^{\prime\prime}$. This gives:
$$g^{\prime}(\rho g^{\prime} + 2m) = 0$$
If $g^{\prime} = 0$, then $g(\rho) = \overline{B}$. Otherwise, $\rho g^{\prime} + 2m = 0$ and so $g(\rho) = \ln(B\rho^{-2m})$ for some constant $B > 0$.

Therefore, when $n\geq 3$, solutions are:
\begin{subequations}
\begin{align}
    \phi(z, \rho) &= (Az^m + B)^{\frac{2}{m}} \, , \; A,B > 0  \\
    \phi(z, \rho) &= (Az^m + B\rho^{-2m})^{\frac{2}{m}} \, , \; A,B > 0  \label{m3sol}
\end{align}
\end{subequations}

For $n = 2$, (\ref{m q eq}) is vacuous and (\ref{m p eq}) gives a linear ODE:
$$\rho g^{\prime\prime} + g^{\prime} = 0$$
So $g(\rho) = \ln(B\vert\rho\vert^C)$, for constants $B > 0$ and $C\in\mathbb{R}$.

Hence, solutions for $n=2$ are:
\begin{equation}
    \phi(z,\rho) =  (Az^m + B\vert\rho\vert^C)^{\frac{2}{m}} \, , \; A,B > 0 \, , \; C\in\mathbb{R}
\end{equation}

For $n = 1$, the Ricci-flat equations are $2(m-1)p_{2m} - q_{2m} = 2(m-1)p_m + q_m = p_0 + q_0 = 0$, by independence of $z$ and $\rho$. As before, since $m > 2$, $2(m-1)p_{2m} - q_{2m} = 0$ implies $f^{\prime} = 0$, and equation $2(m-1)p_m + q_m = 0$ is automatically satisfied. The remaining equation gives:
$$m\rho g^{\prime\prime} - \rho (g^{\prime})^2 - m g^{\prime} = 0$$
If $g^{\prime} = 0$, the above equation is trivially satisfied; then $g(\rho) = \overline{B}$. Otherwise, this is yet again a Bernoulli equation in $g^{\prime}$. Let $u := (g^{\prime})^{-1}$ to obtain a linear ODE:
$$m\rho u^{\prime} + \rho + m u = 0$$
Its solution gives $g(\rho) = \ln(B\vert \rho^2 + C \vert^{-m})$ for constants $B > 0$ and $C\in\mathbb{R}$.

So, for $n=1$, solution are:
\begin{subequations}
\begin{align}
    \phi(z, \rho) &= (Az^m + B)^{\frac{2}{m}} \, , \; A,B > 0  \\
    \phi(z, \rho) &= (Az^m + B\vert \rho^{2} + C \vert^{-m})^{\frac{2}{m}} \, , \; A,B > 0 \, , \; C\in\mathbb{R}
\end{align}
\end{subequations}

For an odd integer $m>2$, all formulas still hold, but the metric generated changes signature according to the sign of $z$, because it determines the sign of $\Lambda$. For $m=2$, $\phi$ simplifies to give a Riemannian metric; in this case, the non-trivial Ricci-flat equations are multiples of the previously found equations for Riemannian metrics.

Finally, taking $\phi(z,\rho) = \left(\e^{f(\rho)}z^{m} - \e^{g(\rho)}\right)^{\frac{2}{m}}$ for some integer $m > 2$ with $m \equiv 2 (\text{mod } 4)$ gives a well-defined metric $L = \alpha^2\phi$ with Lorentz signature $(+,-,\ldots,-)$, since $\Omega = - \frac{2\e^g}{(\e^f z^m - \e^g)^{1-\frac{2}{m}}}$ and $\Lambda = - \frac{4(m-1)\e^{f+g}z^{m-2}}{(\e^f z^m - \e^g)^{2(1-\frac{2}{m})}}$.

In this setting, the Ricci components are:
\begin{align*}
    P &= -\frac{1}{2m^2(m-1)\rho}\left[ p_{2m}\e^{2(f-g)}z^{2m} - p_2\e^{f-g}z^2  + p_0 \right] \\
    Q &= \frac{1}{4m^2(m-1)^2\rho^3}\left[q_{2m}\e^{2(f-g)}z^{2m} + q_m\e^{f-g}z^m - q_0 \right]
\end{align*}
where $p_i$, $q_j$ are as before. Thus, the Ricci-flat equations coincide with the positive-definite case.

With some thought, one might consider these equations for other values of $m$. When $m > 2$ is divisible by $4$, one may take $L = \alpha^m\left(\e^{f(\rho)}z^{m} - \e^{g(\rho)}\right)$ to consider Finsler spacetimes in the sense of Pfeifer and Wohlfarth \cite{pfeifer:wohlfarth}. When $m > 2$ is odd, $F = \alpha\left(\e^{f(\rho)}z^{m} - \e^{g(\rho)}\right)^{\frac{1}{m}}$ already makes sense. However, in both cases, one needs to become concerned with the domain of $z$ and $\rho$ to ensure $\Omega$, $\Lambda$ are defined and their sign give the appropriate signature.
\end{ex}

\begin{ex}[Randers metrics]
Assume $\phi(z,\rho) = (\sqrt{\e^{f(\rho)}z^2 + \e^{g(\rho)}} + \varepsilon\e^{\frac{f(\rho)}{2}}z)^2$ with $0 < \vert \varepsilon \vert < 1$, so $F = \alpha\sqrt{\phi}$ gives a positive-definite Randers metric.

Indeed, $\Omega = 2\left(\frac{\sqrt{\e^f z^2 + \e^g} + \varepsilon\e^{\frac{f}{2}}z}{\sqrt{\e^f z^2 + \e^g}}\right)\e^g$ and $\Lambda = 4\left(\frac{\sqrt{\e^f z^2 + \e^g} + \varepsilon\e^{\frac{f}{2}}z}{\sqrt{\e^f z^2 + \e^g}}\right)^3\e^{f+g}$.

The Ricci curvature components are:
\begin{dgroup*}
\bmatn
P = -\frac{1}{4\rho\sqrt{\e^f z^2+\e^g}(\sqrt{\e^f z^2+\e^g} + \varepsilon\e^{\frac{f}{2}}z)}\left[ p_4\e^{2f-g}z^4 + 2\varepsilon p_3\e^{\frac{3f}{2}-g}\sqrt{\e^f z^2+\e^g}z^3 + p_2\e^f z^2 + \varepsilon p_1\e^{\frac{f}{2}}\sqrt{\e^f z^2+\e^g}z + p_0\e^g \right]
\ematn
\bmatn
Q = \frac{1}{4\rho^3(\e^f z^2+\e^g)^2(\sqrt{\e^f z^2+\e^g} + \varepsilon\e^{\frac{f}{2}}z)^2}\left[ q_6\e^{3f}z^6 + 2\varepsilon q_5\e^{\frac{5f}{2}}\sqrt{\e^f z^2+\e^g}z^5 + 2q_4\e^{2f+g}z^4 + 4\varepsilon q_3e^{\frac{3f}{2}+g}\sqrt{\e^f z^2+\e^g}z^3 + q_2\e^{f+2g}z^2 + 2\varepsilon q_1\e^{\frac{f}{2}+2g}\sqrt{\e^f z^2+\e^g}z + q_0\e^{3g} \right]
\ematn
\end{dgroup*}
where $p_i$, $q_j$ are functions of $\rho$, $f$, $g$ and its derivatives of order up to two. Particularly,
\begin{dgroup*}
\bmatn
p_4 = 2(\varepsilon^2 + 1)\rho f^{\prime\prime} - ((n+1)\varepsilon^2 - 1)\rho(f^{\prime})^2 + (n-2)(\varepsilon^2 + 1)\rho f^{\prime}g^{\prime} + 2(n-1)(\varepsilon^2 + 1)f^{\prime}
\ematn
\bmatn
p_3 = 2\rho f^{\prime\prime} - \frac{1}{4}((n+2)\varepsilon^2 + (n-2))\rho (f^{\prime})^2 + (n-2)\rho f^{\prime}g^{\prime} + 2(n-1)f^{\prime}
\ematn
\end{dgroup*}

For $n\geq 2$, the Ricci-flat equations reduce to $p_i = q_j = 0$, by independence of $z$ and $\rho$. Taking $p_4 - (\varepsilon^2 + 1)p_3 = 0$ reads
$$ \frac{(n+2)}{4}(\varepsilon^2 - 1)^2\rho (f^{\prime})^2 = 0 \, .$$
So $f^{\prime} = 0$, and the remaining equations simplify to:
\begin{equation}
\label{r p eq}
    2\rho g^{\prime\prime} + (n-2)\rho (g^{\prime})^2 + 2(2n-3)g^{\prime} = 0
\end{equation}
\begin{equation}
\label{r q eq}
    (n-2)[2\rho g^{\prime\prime} - \rho(g^{\prime})^2 - 2g^{\prime}] = 0
\end{equation}
Hence, $f(\rho) = \overline{A}$ and $g(\rho)$ must be determined from the above equations.

When $n\geq 3$, one may combine equations (\ref{r p eq}) and (\ref{r q eq}) to eliminate $g^{\prime\prime}$, obtaining:
$$(n-1)(\rho g^{\prime} + 4)g^{\prime} = 0$$
If $g^{\prime} = 0$, then $g(\rho) = \overline{B}$. Otherwise, $\rho g^{\prime} + 4 = 0$ and so $g(\rho) = \ln (B\rho^{-4})$ for some constant $B>0$.

Thus, solutions for $n \geq 3$ are:
\begin{subequations}
\begin{align}
    \phi(z, \rho) &= (\sqrt{A z^2 + B} + \varepsilon\sqrt{A}z)^2 \, , \; A,B > 0 \,  , \; 0 < \vert \varepsilon \vert < 1 \\
    \phi(z, \rho) &= (\sqrt{A z^2 + B\rho^{-4}} + \varepsilon\sqrt{A}z)^2 \, , \; A,B > 0 \, , \; 0 < \vert \varepsilon \vert < 1 \, , \; C\in\mathbb{R} \label{r3sol}
\end{align}
\end{subequations}

For $n = 2$, (\ref{r q eq}) is vacuous and (\ref{r p eq}) becomes a linear ODE:
$$\rho g^{\prime\prime} + g^{\prime} = 0$$
So $g(\rho) = \ln(B\vert\rho\vert^C)$, for constants $B > 0$ and $C\in\mathbb{R}$.

Hence, when $n=2$, solutions are:
\begin{equation}
    \phi(z,\rho) =  \left( \sqrt{A z^2 + B\vert \rho \vert ^{C}} + \varepsilon\sqrt{A}z \right)^2 \, , \; A,B > 0 \, , \; C\in\mathbb{R} , \; 0 < \vert \varepsilon \vert < 1
\end{equation}

Finally, for $n = 1$, the Ricci-flat condition once again implies $f^{\prime} = 0$, although the computation is lengthier and will be omitted. All remaining equations are automatically satisfied.

Therefore, for $n=1$, solutions are:
\begin{equation}
    \phi(z,\rho) = \left( \sqrt{A z^2 + \e^{g(\rho)}} + \varepsilon\sqrt{A}z \right)^2 \, , \; A > 0 \, , \; g\in C^{\infty} \, , \; 0 < \vert \varepsilon \vert < 1
\end{equation}

Clearly, one may rewrite solutions as $\phi(z,\rho) = \left( \sqrt{A z^2 + \e^{g(\rho)}} + D z \right)^2$ for any constant $D$ satisfying $D^2A^{-1} < 1$. More generally, it is possible to look for solutions in the form $\phi(z,\rho) = \left( \sqrt{\e^{f(\rho)} z^2 + \e^{g(\rho)}} \pm h(\rho) z \right)^2$ with $h^2(\rho) < \e^{f(\rho)}$, but the calculations quickly become cumbersome. In either case, it is uncertain how to consider Lorentz signature (if possible).
\end{ex}

\section{Discussion}

The Hessian of the Ricci curvature
$$\ric_{AB} = \frac{1}{2}\left[\ric \right]_{y^A y^B}$$
was the first notion for Ricci curvature tensor of Finsler metrics introduced by Akbar-Zadeh in 1988. Evidently, $\ric_{AB} = 0$ if and only if $\ric = 0$, and they imply the vanishing of the scalar curvature $R = g^{AB}\ric_{AB}$. By defining the modified Einstein tensor
$$ G_{AB} = \ric_{AB} - \frac{1}{2}g_{AB}R $$
in \cite{li:chang}, Li an Chang established the equivalence between the vacuum field equation for Finsler spacetime and the vanishing of the Ricci curvature. However, the notion of Ricci curvature tensor for Finsler metrics is not unique. If $R^{\;A}_{B\;CD}$ is the Riemann curvature tensor for Finsler metrics introduced by Berwald in 1926, then
$$\widetilde{\ric}_{AB} = \frac{1}{2}\left( R^{\;C}_{A\;CB} + R^{\;C}_{B\;CA} \right)$$
is another notion of Ricci curvature tensor introduced by Li and S. in \cite{li:shen}. Moreover, these Ricci tensors differ by a non-Riemannian quantity; namely,
$$ \widetilde{\ric}_{AB} - \ric_{AB} = H_{AB} = \frac{1}{2}\left( [\chi_B]_{y^A} + [\chi_A]_{y^B} \right) \, , $$
where the $\chi$-curvature tensor is given by
$$\chi_A = \frac{1}{2}\left[ \Pi_{x^B y^A}y^B - \Pi_{x^A} - 2\Pi_{y^A y^B}G^B \right]$$
with $\Pi = \frac{\partial G^C }{\partial y^C}$. So $\widetilde{\ric}_{AB} = 0$ if and only if $\ric_{AB} = 0$ and $H_{AB} = 0$; in words, the vanishing of $\widetilde{\ric}_{AB}$ is a stronger condition than the vanishing of $\ric_{AB}$. In particular, if $\ric = 0$ and $\chi_A = 0$, then $\widetilde{\ric}_{AB} = 0$.

For the proposed metrics $F = \alpha\sqrt{\phi(z,\rho)}$, we have $\Pi = \Psi(x^m y^m)$, where
\begin{equation}
\label{Psi}
    \Psi := U_z + (n+2)V + (n-1)W \, ,
\end{equation}
and the $\chi$-curvature is
\begin{dgroup}
\bmat
\chi_0 = \left[ \frac{1}{2\rho}\Psi_{z\rho} - U\Psi_{zz} - W\Psi_{z} \right]\frac{(x^my^m)^2}{\alpha} + \frac{1}{2}(2\rho^2 W + 1)\Psi_z\alpha
\emat
\bmat
\chi_i = \left[ zU\Psi_{zz} - \frac{z}{2\rho}\Psi_{z\rho} + (U + 2zW)\Psi_z \right] \frac{(x^m y^m)^2}{\alpha^2}y^i -\frac{z}{2}(2\rho^2 W + 1)\Psi_z y^i - (U + zW)\Psi_z(x^m y^m)x^i
\emat
\end{dgroup}
Clearly, $\Psi_z = 0$ is a sufficient condition for the vanishing of the $\chi$-curvature. By direct verification, all solutions in previous section satisfy $\Psi_z = 0$. Thus, they are strongly Ricci-flat metrics: $\widetilde{\ric}_{AB} = \ric_{AB} = 0$.

In addition to the examples presented here, it seems to be feasible (although lengthy) to construct other types of (strongly) Ricci-flat metrics in the proposed form; particularly, one may look for series expansions. The same type of construction also seems to work well for Ricci-isotropic metrics, $\ric = [(n+1)-1] k(x)F^2$. At the very least the PDE characterization is similar to describe; namely, for $n \geq 2$, $F=\alpha\sqrt{\phi(z,\rho)}$ is Ricci isotropic if and only if $P = n k \phi$ and $Q=0$. It might be wise, however, to spend such efforts with a wider class of warped product Finsler metrics, which may allow for global solutions on $\mathbb{R}\times M$; for instance, a class of Finsler metrics defined by
\begin{equation}
    F = \alpha \sqrt{\phi(z,\overline{x})} \, ,
\end{equation}
for $\alpha$ any Riemannian metric on $M$, $z$ as before and $\phi$ some appropriate function on $\mathbb{R}\times M$.

\appendix
\section{Derivatives}

Derivatives of $F^2$:
\begin{dgroup*}
\bmatn
    \left[ F^2 \right]_{y^0} = \alpha\phi_z
\ematn
\bmatn
    \left[ F^2 \right]_{y^i} = \Omega y^i
\ematn
\bmatn
    \left[ F^2 \right]_{x^0} = 0
\ematn
\bmatn
    \left[ F^2 \right]_{x^i} = \frac{1}{\rho}\phi_{\rho}\alpha^2 x^i
\ematn
\bmatn
    \left[ F^2 \right]_{y^A x^0} = 0
\ematn
\bmatn
    \left[ F^2 \right]_{y^0 x^i} = \frac{1}{\rho}\phi_{z\rho}\alpha x^i
\ematn
\bmatn
    \left[ F^2 \right]_{y^i x^j} = \frac{1}{\rho}\Omega_{\rho}x^j y^i
\ematn
\end{dgroup*}

Derivatives of $G^A$:
\begin{dgroup*}
\bmatn
    \left[ G^A \right]_{x^0} = 0
\ematn
\bmatn
    \left[ G^0 \right]_{x^i} = \frac{1}{\rho}(U_{\rho} + z V_{\rho})x^i(x^m y^m)\alpha + (U + z V)y^i\alpha
\ematn
\bmatn
    \left[ G^j \right]_{x^i} = \frac{1}{\rho}(V_{\rho} + W_{\rho})x^i y^j(x^m y^m) + (V + W)y^i y^j - \frac{1}{\rho}W_{\rho}x^i x^j\alpha^2 - W\delta_i^j\alpha^2
\ematn
\bmatn
    \left[ G^0 \right]_{y^0} = (U_z + V + z V_z)(x^m y^m)
\ematn
\bmatn
    \left[ G^j \right]_{y^0} = (V_z + W_z)(x^m y^m)\frac{y^j}{\alpha} - W_z x^j\alpha
\ematn
\bmatn
    \left[ G^0 \right]_{y^i} = (U -z U_z - z^2 V_z)(x^m y^m)\frac{y^i}{\alpha} + (U + z V)x^i\alpha
\ematn
\bmatn
    \left[ G^j \right]_{y^i} = (V + W)(x^m y^m)\delta_i^j - z(V_z + W_z)(x^m y^m)\frac{y^i y^j}{\alpha^2} + (V + W)x^i y^j + (z W_z - 2W)x^j y^i
\ematn
\bmatn
    \left[ G^B \right]_{x^0 y^A} = 0
\ematn
\bmatn
    \left[ G^0 \right]_{x^i y^0} = \frac{1}{\rho}(U_{z\rho} + V_{\rho} + z V_{z\rho})(x^m y^m)x^i + (U_z + V + z V_z)y^i
\ematn
\bmatn
    \left[ G^j \right]_{x^i y^0} = \frac{1}{\rho}(V_{z\rho} + W_{z\rho})(x^m y^m)\frac{x^i y^j}{\alpha} + (V_z + W_z)\frac{y^i y^j}{\alpha} - \frac{1}{\rho}W_{z\rho}x^i x^j\alpha - W_z\delta_i^j\alpha
\ematn
\bmatn
    \left[ G^0 \right]_{x^i y^j} = (U + z V)\delta_j^i\alpha + (U - z U_z - z^2 V_z)\frac{y^i y^j}{\alpha} + \frac{1}{\rho}(U_{\rho} - z U_{z\rho} - z^2 V_{z\rho})(x^m y^m)\frac{x^i y^j}{\alpha} + \frac{1}{\rho}(U_{\rho} + z V_{\rho})x^i x^j\alpha
\ematn
\bmatn
    \left[ G^k \right]_{x^i y^j} = (V + W)(\delta_j^i y^k + \delta_j^k y^i) + (z W_z - 2W)\delta_i^k y^j + \frac{1}{\rho}(V_{\rho} + W_{\rho})(x^m y^m)x^i\delta_j^k + \frac{1}{\rho}(V_{\rho} + W_{\rho})x^i x^j y^k + \frac{1}{\rho}(z W_{z\rho} - 2W_{\rho})x^i x^k y^j -\frac{z}{\rho}(V_{z\rho} + W_{z\rho})(x^m y^m)\frac{x^i y^j y^k}{\alpha^2} - z(V_z + W_z)\frac{y^i y^j y^k}{\alpha^2}
\ematn
\bmatn
    \left[ G^0 \right]_{y^0 y^0} = (U_{zz} + 2 V_z + z V_{zz})\frac{(x^m y^m)}{\alpha}
\ematn
\bmatn
    \left[ G^k \right]_{y^0 y^0} = (V_{zz} + W_{zz})(x^m y^m)\frac{y^k}{\alpha^2} - W_{zz}x^k
\ematn
\bmatn
    \left[ G^0 \right]_{y^0 y^i} = (U_z + V + z V_z)x^i - z(U_{zz}+ 2\psi_z + z V_{zz})(x^m y^m)\frac{y^i}{\alpha^2}
\ematn
\bmatn
    \left[ G^k \right]_{y^0 y^i} = (V_z + W_z)\frac{(x^m y^m)}{\alpha}\delta_i^k + (V_z + W_z)\frac{x^i y^k}{\alpha} + (z W_{zz} - W_z)\frac{x^k y^i}{\alpha} - (V_z + z V_{zz} + W_z + z W_{zz})\frac{(x^m y^m)}{\alpha}\frac{y^i y^k}{\alpha^2}
\ematn
\bmatn
    \left[ G^0 \right]_{y^i y^j} = (U - z U_z - z^2 V_z)\frac{(x^m y^m)}{\alpha}\delta_j^i + (U - z U_z - z^2 V_z)\left[ \frac{x^i y^j}{\alpha} + \frac{x^j y^i}{\alpha}\right] + (-U + z U_z + z^2 U_{zz} + 3z^2 V_z + z^3 V_{zz})\frac{(x^m y^m)}{\alpha}\frac{y^i y^j}{\alpha^2}
\ematn
\bmatn
    \left[ G^k \right]_{y^i y^j} = z(3(V_z + W_z) + z(V_{zz} + W_{zz}))(x^m y^m)\frac{y^i y^j y^k}{\alpha^4} - z(V_z + W_z)\left[ \frac{x^i y^j y^k}{\alpha^2} + \frac{x^j y^i y^k}{\alpha^2} \right] - z(z W_{zz} - W_z)\frac{x^k y^i y^j}{\alpha^2} + (V + W)(x^i \delta^j_k + x^j \delta^k_i) + (z W_z - 2W)x^k\delta^i_j - z(V_z + W_z)(x^m y^m)\left[ \frac{\delta^j_i y^k}{\alpha^2} + \frac{y^i\delta^j_k}{\alpha^2} + \frac{\delta^k_i y^j}{\alpha^2}\right]
\ematn
\end{dgroup*}

Derivatives of $\Pi$:
\begin{dgroup*}
\bmatn
\Pi_{x^0} = 0
\ematn
\bmatn
\Pi_{x^i} = \Psi y^i + \frac{1}{\rho}\Psi_{\rho}(x^m y^m)x^i
\ematn
\bmatn
\Pi_{y^0} = \Psi_{z}\frac{(x^m y^m)}{\alpha}
\ematn
\bmatn
\Pi_{y^j} = -z\Psi_z\frac{(x^m y^m)}{\alpha}\frac{y^j}{\alpha} + \Psi x^j
\ematn
\bmatn
\Pi_{x^0 y^A} = 0
\ematn
\bmatn
\Pi_{x^i y^0} = \Psi_z\frac{y^i}{\alpha} + \frac{1}{\rho}\Psi_{z\rho}\frac{(x^m y^m)}{\alpha}x^i
\ematn
\bmatn
\Pi_{x^i y^j} = \Psi\delta_i^j - z\Psi_z\frac{y^i y^j}{\alpha^2} - \frac{z}{\rho}\Psi_{z\rho}\frac{(x^m y^m)}{\alpha}x^i\frac{y^j}{\alpha} + \frac{1}{\rho}\Psi_{\rho}x^ix^j
\ematn
\bmatn
\Pi_{y^0 y^0} = \Psi_{zz}\frac{(x^m y^m)}{\alpha^2}
\ematn
\bmatn
\Pi_{y^0 y^j} = \Psi_z\frac{x^j}{\alpha} -\left(\Psi_z + z\Psi_{zz} \right)\frac{(x^m y^m)}{\alpha^2}\frac{y^j}{\alpha}
\ematn
\bmatn
\Pi_{y^i y^j} = -z\Psi_z\left[ \frac{(x^m y^m)}{\alpha^2}\delta_{i}^{j} + \frac{x^i y^j}{\alpha^2} + \frac{x^j y^i}{\alpha^2} \right] + z\left( 3\Psi_z + z\Psi_{zz} \right)\frac{(x^m y^m)}{\alpha^2}\frac{y^i y^j}{\alpha^2}
\ematn
\end{dgroup*}

\section*{Acknowledgement}
This project is supported in part by NNSFC (No. 11671352) and CNPq (No. 217974/2014-7).

\bibliography{mybib}

\end{document}